\documentclass[12pt,oneside]{article}
\usepackage{amsmath,amssymb,amsfonts,amsthm}
\usepackage{amscd}
\usepackage{graphicx}
\usepackage[all]{xy}
\usepackage{multirow}
\usepackage{rotating}
\usepackage{hyperref}
\usepackage{hyperref}
\hypersetup{
    colorlinks=true,
    linkcolor=blue,
    filecolor=magenta,
    urlcolor=cyan,
    citecolor=blue}
\usepackage{lscape}
\newtheorem{thm}{Theorem}[section]
\newtheorem{cor}[thm]{Corollary}
\newtheorem{lem}[thm]{Lemma}
\newtheorem{prop}[thm]{Proposition}
\newtheorem{defn}[thm]{Definition}
\newtheorem{example}{Example}
\newtheorem{rem}[thm]{Remark}

\newcommand{\T}{{\cal T}}

\newcommand{\Real}{\mathbb R}

\newcommand{\set}[1]{\left\{#1\right\}}

\newcommand{\To}{\longrightarrow}

\numberwithin{equation}{section}

\begin{document}
\title{ A Concurrent Generalized Kropina Change }
\author{A. Soleiman$^{1}$ and Ebtsam H. Taha$^2$ }
\date{}
\maketitle  
\begin{center}
{$^{1}$ Department of Mathematics, College of Science and Arts - Qurayyat,
Al Jouf University, Skaka,  Kingdom of Saudia Arabia\\
\vspace{0.2cm}
E-mails: amr.hassan@fsc.bu.edu.eg, amrsoleiman@yahoo.com}
\end{center}  
\begin{center}
$^{2}$ Department of Mathematics, Faculty of Science, Cairo University, Giza, Egypt,
\end{center}
\vspace{-0.8cm}
\begin{center}
E-mails: ebtsam.taha@sci.cu.edu.eg, ebtsam.h.taha@hotmail.com
\end{center} 
 
\begin{center}
\textbf{\textit{To the memory of Professor Nabil L. Youssef}}
\end{center}

\maketitle 

\noindent{\bf Abstract.} This paper investigates a generalized Kropina metric featuring a specific $\pi$-form. Start with a Finsler manifold $(M,F)$ admits a concurrent $\pi$-vector field $\overline{\varphi}$, then, examine the $\phi$-concurrent generalized Kropina change defined by $\widehat{F}=\frac{F^{m+1}}{\Phi^{m}}$, where $\Phi$ represents the corresponding $1$-form,  $\Phi^{m}>0$.  We investigate the fundamental geometric objects associated with $\widehat{F}$ in an intrinsic manner after adopting this modification and present an example of a Finsler metric that admits a concurrent vector field along with $\widehat{F}$.   Also, we prove that the geodesic sprays of $F$ and $\widehat{F}$ can never be projectively related.  Moreover,  we show $\overline{\varphi}$ is not concurrent with respect to $\widehat{F}$.  Eventhough, we give a sufficient condition for $\overline{\varphi}$ to be   concurrent with respect to $\widehat{F}$. Finally,  we prove that the  $\phi$-concurrent generalized Kropina change ($F  \longrightarrow \widehat{F}$) preserves the almost rational property of the initial Finsler metric ${F}$. 

\bigskip
\medskip\noindent{\bf Keywords:\/}\, Generalized Kropina metric, Concurrent vector field,  almost rational Finsler metric, Berwald connection

\medskip
\noindent{\bf MSC 2020}: 53C60, 53B40, 58B20.

\vspace{30truept}\centerline{\Large\bf{Introduction}}\vspace{12pt}
The  Kropina metric $F= \frac{\alpha^{m+1}}{\beta^m}; \,\, m=1$ is an interesting metric in Finsler geometry which has been investigated  firstly in \cite{Kropina}.  Also,  for values of $m \neq 0,1$ the Finsler metric  $F= \frac{\alpha^{m+1}}{\beta^m}$ is called generalized Kropina metric which is considered as a special $(\alpha , \beta)$-metrics. The generalized Kropina change of a Finsler metric $F$ (which is not necessarily Riemannian) $F  \longrightarrow \frac{F^{m+1}}{\beta^m}$  has been done for certain Finsler metrics such as,  the m-th root metrics and exponential $(\alpha , \beta)$-metric in  \cite{gen. Krop. change,  gen. Krop. change2017, Tiwari-Kumar-Tayebi, expKrpoina2024}.  Nice results are obtained, for example, the conditions under which the  Kropina change of generalized m-th root metrics are locally projectively flat and locally dually flat. This motivates us to study the generalized Kropina change of an arbitrary Finsler metric which is the main topic of this paper.  This transformation allows us to  investigate the  behaviours of this change on various geometric objects. This change is maybe useful in the context of Finslerian modification of general theory of relatively,  since the generalized Kropina metric has  been used  effectively in \cite{Kropina-Eur. Phys. J. C}.

Special Finsler spaces investigated globally in \cite{amr3, amr2, Hsca, Numata1, r99} and locally  in \cite{Kropina-Eur. Phys. J. C, arFinslerTaha, Taha_Gen_m_Kropina, Tamim-Youssef}.  Both local and global point of views are useful in treating problems in Finsler geometry.  In this paper,  we examine a Finsler manifold $(M,F)$ that possesses a concurrent $\pi$-vector field $\overline{\varphi}$.  We use its corresponding $\pi$-form $\phi:=i_{\overline{\varphi}} g$, where $g$ represents the metric tensor of $F$, leading to the associated function $\Phi(x,y):=\phi(\overline{\eta})$. Next, we examine what we called the $\phi$-concurrent generalized Kropina change $\widehat{F}=F^{m+1}\Phi^{-m}$.  In this context, we compute intrinsic geometric objects related to $\widehat{F}$.  Namely, the supporting form $\widehat{\ell}$, the angular metric tensor $\widehat{\hbar}$, the Finsler metric $\widehat{g}$,  the Cartan torsion $\widehat{\mathbf{T}}$, the geodesic spray $\widehat{G}$, the nonlinear connection $\widehat{\Gamma}$, the Berwald connection and the curvature tensor $\widehat{\Re}$ associated with  $\widehat{\Gamma}$ are identified in terms of the corresponding geometric objects of $F$.  Furthermore, we characterise the non-degenerate property of the metric tensor $\widehat{g}$  in \S 2. 

\par We have noted that the above mentioned geometric objects are not invariant under the $\phi$-concurrent generalized Kropina change except the vertical counterpart of Berwald connection. We find a sufficient condition that makes some geometric objects, namely,  $\widehat{\Gamma} ,\, \Re$ and the horizontal counterpart of Berwald connection to be invariant.  Consequently, we prove that the $\pi$-vector field $\overline{\varphi}$ is not concurrent with respect to $\widehat{F}$,  however, it will be under certain condition. 
On the other hand, we conclude the geodesic sprays $G$ and  $\widehat{G}$ can never be projectively related.   An example which represent our change is provided. We end this work by study the effect of the $\phi$-concurrent generalized Kropina change on an almost rational Finsler metric and obtain interesting results.

\section{Preliminaries}
\par 
Let $M$ be an $n$-dimensional smooth manifold and  $\pi: T M\longrightarrow M$ its  tangent bundle. The vertical subbundle $V(TM)$ is defined to be $ \ker (d\pi)$.  We denote the pullback bundle of the tangent bundle by $\pi^{-1}(T M)$.  Further,  $\mathfrak{F}(TM)$ denotes the algebra of smooth functions on $TM$ and $\mathfrak{X}(\pi )$ the $\mathfrak{F}(TM)$-module of differentiable sections of $\pi^{-1}(T M)$. The elements of $\mathfrak{X}(\pi )$ will be called $\pi$-vector fields and denoted by barred letters $\overline{X}$. \\
 \vspace*{-0.4cm}
\par
We have the short exact sequence \cite{r21, r94a}
$$0\longrightarrow
 \pi^{-1}(TM)\stackrel{\gamma}\longrightarrow TTM \stackrel{\rho}\longrightarrow
\pi^{-1}(TM)\longrightarrow 0 ,\vspace{-0.1cm}$$ where $\gamma$ is the natural injection (which is an isomorphism from $ \pi^{-1}(TM)$ to $V(TM)$) and $\rho := (\pi_{TM}, d \pi)$. The tangent structure $J$  is $(1,1)$-type tensor $J: TTM \longrightarrow TTM$ defined by $J=\gamma \circ \rho$. For all  $ f \in \mathfrak{F}(TM), \, W \in \mathfrak{X}(TM)$, $J$ satisfies:
 \begin{equation}\label{fJformula}
[f W, J]=f[W,J]+df \wedge i_{W}J-d_{J}f\otimes W,
\end{equation}
\begin{equation}\label{J properties}
 i_{ \mathcal{C}} \,J=0 \quad \text{    and    } \,\, \, [\mathcal{C}, J] = -J,\qquad \mathcal{C}:=\gamma\, \overline{\eta},
 \end{equation}
where $\overline{\eta}(u)=(u,u)$ for all $u$ in the slit tangent bundle $\T M:=TM/\set{0}.$ The vector field $\mathcal{C}$ is called Liouville vector filed.
\par
For a linear connection  $D$ on $\pi^{-1}(TM)$, we have    $K:TTM \longrightarrow
\pi^{-1}(TM)$ which is defined by $K(W) =D_W \overline{\eta}$.  Thereby,  the horizontal space at $u \in TM$ is $H_u (TM):= \{ W \in T_u
(TM) \,|\, K(W)=0 \}$.  The connection $D$ is said to be regular if for all $u\in TM$, we have $ T_u (TM)=V_u (TM)\oplus H_u (TM).$ For a regular connection $D$,  the vector bundle
   maps $  \rho |_{H(TM)}$ and $K |_{V(TM)}$
 are isomorphisms. In this case, the map  $\beta:={\rho^{-1}} |_{H(TM)}$
 is called the horizontal map of $D$. A well-Known regular connection is Berwald connection  \cite{r21} which can be defined by \cite[Proposition 4.4]{r92} 
\begin{equation}\label{Berwald} \gamma {D^\circ}_{\beta \rho{Z}} \, \overline{W}:=\gamma \, K[\beta \rho Z,JW], \qquad {D^\circ}_{\gamma\overline{Z}} \, \rho {W}:=\rho[\gamma \overline{Z},{\beta} \overline{W}].
\end{equation}   
\par Moreover,  a spray on $M$ is a smooth vector field $G$ on $\T M$ such that  $J G = \mathcal{C}$ and $[\mathcal{C},G]= G$.  It is clear that $G=\beta\, \overline{\eta}$. A nonlinear connection on $M$ is a vector $1$-form $\Gamma$ on $TM$ which is
smooth on $\T M$ and continuous  on $TM$ such that $J \Gamma=J \text{ and } \,\, \Gamma J=-J $ \cite{r21}.  Consequently, the horizontal and vertical projectors
 associated with $\Gamma$ are
given,  respectively, by
 \begin{equation}\label{hor. and ver. proj.}
h:=\frac{1}{2} (I+\Gamma) = \beta \, \rho,\quad  \quad v:=\frac{1}{2}
 (I-\Gamma)=\gamma \, K.
\end{equation}
Consequently, we get $vJ=J $ and $Jv=0$.
The curvature of $\Gamma$ is defined by
 $\mathfrak{R}:=-\frac{1}{2}[h,h]$, which can be computed 
 using Fr\"{o}licher-Nijenhuis  bracket $[\mathbb{K},\mathbb{L}]$ of two vector 1-forms $\mathbb{K}$ and $\mathbb{L}$ as follows  \cite{r20}:
\begin{eqnarray}\label{brackect of two v.forms}
  [\mathbb{K},\mathbb{L}](W,Z)&=& [\mathbb{K}W,\mathbb{L}Z]+[\mathbb{L}W,\mathbb{K}Z]+\mathbb{K} \mathbb{L}[W,Z]+\mathbb{L} \mathbb{K}[W,Z] \nonumber \\ 
   && -\mathbb{K}[\mathbb{L} W,Z]-\mathbb{K}[ W,\mathbb{L} Z]-\mathbb{L}[\mathbb{K} W,Z]-\mathbb{L}[ W,\mathbb{K} Z].
\end{eqnarray}
In particular,   the Nijenhuis torsion $N_{\mathbb{L}}$ of a vector $1$-form  $\mathbb{L}$ is  defined by
\begin{equation}\label{Nk}
     N_{\mathbb{L}}(W,Z):= \frac{1}{2}[\mathbb{L},\mathbb{L}](W,Z)=[\mathbb{L}W,\mathbb{L}Z]+\mathbb{L}^{2}[W,Z] -\mathbb{L}[\mathbb{L} W,Z]-\mathbb{L}[ W,\mathbb{L} Z].
\end{equation}
It is known that $N_{J} = 0$ and $J^{2} =0$ which imply 
\begin{equation}\label{JJ}
[JW,JZ]=J[W,JZ]+J[JW,Z].
\end{equation}
\begin{defn}
A Finsler metric on $M$ is a function 
$F: TM \To [0,\infty) $  such~that{\em:}
 \begin{description}
    \item[(a)] $F $ is  $C^{\infty}$ on  $\T M$ and  $C^{0}$ on $TM$,
    \item[(b)]$F$ satisfies  $\mathcal{L}_{\mathcal{C}} F=F$, where $\mathcal{L}_\mathcal{C}$ is the Lie derivative in the direction of $\mathcal{C}$,
\item[(c)] the Hilbert $2$-form
    $dd_{J}E$  has a maximal rank, where $E=\frac{1}{2}F^{2} $ is the Finsler energy function. 
 \end{description}
The Finsler metric tensor $g$ induced by $\,F\,$ on $\pi^{-1}(TM)$  is defined as follows\emph{ \cite{r94a}}
\begin{equation}\label{g}
g(\rho W,\rho Z):=dd_{J}E(JW,Z), \ \forall W,Z \in
    \mathfrak{X}(TM).
\end{equation}
  In this case, the pair $(M,F)$ is called  a Finsler manifold and $F$ is  a Finsler metric.\\
 If the above conditions are satisfied on a conic subbundle $\mathbf{D}$ of $TM$,  then  $(M,F)$ is called a conic pseudo-Finsler manifold.  A conic subbundle of $TM$ is a non-empty open subset $\mathbf{D} \subset \T M$ that is invariant under the scaling of its tangent vectors by positive real numbers.  Additionally, it must satisfy the condition $\pi(\mathbf{D}) = M $. Furthermore, it is assumed that for every  $x \in M$, the set $ \mathbf{D}_{x} := \mathbf{D} \cap T_{x}M$  is a connected set.
 \end{defn}
 One can easily note that,  $g(\overline{W},\overline{Z})=dd_{J}E(\gamma \overline{W}, \beta \overline{Z})$ for all $\overline{Z}, \overline{W} \in \mathfrak{X}(\pi ).$
 The normalized supporting element $\ell$ (or supporting form) is   defined by $\ell =F^{-1}i_{\overline{\eta}}\:g$ and the angular metric tensor $\hbar := g-\ell \otimes \ell$.  Moreover, the geodesic spray $G$ of $F$ satisfies
    $i_{G}\,dd_JE =-d E.$ Additionally, the Barthel connection $\Gamma$ can be written in terms of $G$ as $\Gamma = [J,G]$  \cite{r21}.  Also, the Berwald connection satisfies  \cite{r92}
\begin{equation}\label{metricity}
(D^{\circ}_{\beta \overline{X}}g)(\overline{Y},\overline{Z})=-2\widehat{\mathbf{P}}(\overline{X},\overline{Y},\overline{Z}),\quad (D^{\circ}_{\gamma \overline{X}}\,g)(\overline{Y},\overline{Z})=2\mathbf{T}(\overline{X},\overline{Y},\overline{Z}).
\end{equation}

\section{ $\phi$-concurrent generalized  Kropina change}
This section starts with  an intrinsic examination of what we call the $\phi$-concurrent generalized Kropina change $F \longrightarrow \widehat{F}$. This study investigates the relationship between the supporting forms ($\ell$ and $\widehat{\ell}$), the angular metric tensors ($\hbar$ and $\widehat{\hbar}$), Finsler metric tensors ($g$ and $\widehat{g}$) and the Cartan torsions ($\mathbf{T}$ and $\widehat{\mathbf{T}}$) associated with this transformation. Additionally,  we provide the condition that makes $\widehat{F}$ non-degenerate. Then, we give an example of our $\phi$- concurrent generalized Kropina change. In addition, we continue to examine the expression of its geodesic spray $\widehat{G}$  in relation to the geodesic spray $G$ of $F$ and prove that $G$ and $\widehat{G}$ cannot be projectively related. Furthermore, the relation between the two Barthel connections ($\Gamma$ and $\widehat{\Gamma}$) is established, along with the derivation of the relations between the Barthel curvature tensors ($\Re$ and $\widehat{\Re}$), as well as the Berwald connections ($D^\circ$ and $\widehat{D^\circ}$).  We conclude that the $\pi$-vector field $\overline{\varphi}$ is not concurrent for $\widehat{F}$ (in general) but it will be under certain condition.
\begin{defn}\emph{\cite{r94a}} \label{concurrent} Given a Finsler manifold $(M,F)$. A non-vanishing $\pi$-vector field  $\overline{\varphi}$ is
said to be  a concurrent $\pi$-vector field if 
\begin{equation}\label{ch512.eq.1}
  {{D}}^{\circ}_{\beta \overline{W}}\,\overline{\varphi} =- \overline{W} , \quad \quad
    {{D}}^{\circ}_{\gamma \overline{W}}\,\overline{\varphi}=0.
 \end{equation}
 \end{defn}
Therefore,  it associated $\pi$-form $\phi:=i_{\overline{\varphi}}\,g$   satisfies
\begin{equation}
({{D}}^{\circ}_{\beta  \overline{W}}\phi)(\overline{Z}) =-g(\overline{W},\overline{Z}), \quad  ({{D}}^{\circ}_{\gamma \overline{W}}\phi)(\overline{Z}) =0.
\end{equation}
\begin{rem}
\par Let us fix our notation throughout the entire paper:\\
$\bullet$  $\overline{\varphi}$ denotes a concurrent $\pi$-vector field with respect to $F$,  \\
$\bullet$  $\phi$ is the $\pi$-form  associated with $\overline{\varphi}$,  \\
$\bullet$ $\Phi $ is a smooth function on TM corresponding to $\phi$ defined at each point by $${\Phi }:=g(\overline{\varphi},\overline{\eta})=\phi(\overline{\eta}).$$
$\bullet$  $||\overline{\varphi}||_{g}:=g(\overline{\varphi},\overline{\varphi})=\phi(\overline{\varphi})$ is the length of $ \overline{\varphi}$ with respect to $F$.
\end{rem}
\begin{lem}\label{independent of the directional argument}
A concurrent $\pi$-vector field $\overline{\varphi}$ on a Finsler manifold $(M,F)$ and its corresponding $\pi$-form $\phi$ have no  dependence of  the directional argument $y$ \emph{ \cite[Theorem 3.7]{r94a}}.   That is,  $$D^{\circ}_{\gamma \overline{W}}\overline{\varphi}=0 =D^{\circ}_{\gamma \overline{W}}\,  \phi,  \quad \forall\, \overline{W} \in \mathfrak{X}(\pi ).$$  Consequently, we posses 
\begin{equation}\label{P properties}
 i_{\gamma \overline{\varphi}} \,J  =0, \qquad d_{J} ||\overline{\varphi}||_{g}=0, \qquad [\gamma \overline{\varphi}, J]=0.
\end{equation}
\end{lem}
\begin{defn}  Given a Finsler manifold $(M,F)$ equipped with a concurrent $\pi$-vector field $\overline{\varphi}$ with the corresponding $\pi$-form $\phi$ and the function $\Phi =\phi(\overline{\eta})$.  Define  \begin{equation}\label{change}
\widehat{F}=F^{m+1}\Phi^{-m},
\end{equation}
with $m \neq 0, -1$.  The change ($F \longrightarrow \widehat{F}$) is called the $\phi$-concurrent generalized Kropina change. If $\widehat{F}$ is a conic pseudo-Finsler metric on $M$,  then it is going to be called  a $\phi$-concurrent generalized Kropina metric.
\end{defn}                                       
\begin{lem}\label{B}  \emph{\cite{r94a, square metric, Soleiman-Taha_Mat}} Consider a Finsler manifold which admits a concurrent $\pi$-vector field $\overline{\varphi}$.  For all $X \in  \mathfrak{X}(TM)$ and $\overline{Z}, \overline{W} \in \mathfrak{X}(\pi )$, the following hold:
\begin{description}
   \item[(a)]$d_{J}\Phi(\gamma \overline{W})= 0, \quad  d_{J}\Phi(\beta \overline{W})={{D}}^{\circ}_{\gamma \overline{W}}\Phi = \phi(\overline{W}),\, \,\,d\Phi(X)=\phi(K X)-F\ell({\rho X})$,
        \item[(b)]$d_{J}\,F(\gamma \overline{W})= 0,  \,\, d_{J}F(\beta \overline{W})={{D}}^{\circ}_{\gamma \overline{W}}\,F= \ell(\overline{W}),\,\,dF(X)=dF(\gamma K X)=\ell(K X),$
             \item[(c)]$d_{h}\,\Phi(\beta \overline{W})=d \Phi(\beta \overline{W})= {{D}}^{\circ}_{\beta \overline{W}}\,\Phi=-F\,\ell(\overline{W}) $.\\ In particular,    $d \Phi(G)=-F^2, \qquad (D^{\circ}_{G} \,\phi)(\overline{W})= -g(\overline{W},\overline{\eta})=-F\,\ell(\overline{W})$,
        \item[(d)]$d_{h}\,F(\beta \overline{W})= d F(\beta \overline{W})={{D}}^{\circ}_{\beta \overline{W}}\,F=0 $,
       \item[(e)] $({{D}}^{\circ}_{\gamma \overline{W}} \,\ell)(\overline{Z})=F^{-1} \hbar(\overline{W},\overline{Z}),\quad \,(D^{\circ}_{G}\, \ell)(\overline{W})= 0,  \quad \rho[G,X]=D^{\circ}_G \rho X-K X$,
         \item[(f)] a smooth function $f$ of two variables $ F \text{ and }\Phi$ satisfies
   \begin{equation}\label{f,b}
   {D^\circ_{\gamma \overline{W}}}f(F,\Phi)=d_Jf(\beta \overline{W})=\frac{\partial f}{\partial F} \, \ell(\overline{W})
        +\frac{\partial f}{\partial \Phi} \, \phi(\overline{W}).
   \end{equation}
     \end{description}
  \end{lem}
  \begin{prop}\label{hh1}Consider a Finsler manifold $(M,F)$  equipped with a concurrent $\pi$-vector field  $\overline{\varphi}$.  Under the $\phi$-concurrent generalized Kropina  change {\em (\ref{change})}, we obtain:  \,
\begin{description}
\item[(1)] The vertical counterpart for Berwald connection
 ${{D}}^{\circ}_{\gamma \overline{X}} \overline{Y}$ is invariant, i.e.,
 \begin{equation}\label{D}
 {{\widehat{D}}}^{\circ}_{\gamma \overline{X}}\, \overline{Y}= {{D}}^{\circ}_{\gamma \overline{X}}\, \overline{Y}.
 \end{equation}
\item[(2)] The total derivative of the Finsler energy functions $d\widehat{E}$ and $dE$ are related by
\begin{eqnarray}\label{AB}
   d\widehat{E}
 &=& \Phi^{-2 m-1} F^{2 m+1}\,\set{(m+1)\,\Phi\,dF-m \,F\, d\Phi }.\end{eqnarray}
 
  \item[(3)]  The supporting elements $\widehat{\ell}$ and $\ell$ are related by
  \begin{equation}\label{ell}
     \widehat{\ell}(\overline{X})={F}^{m}{\Phi}^{-m}\set{(m+1)\ell(\overline{X})
     -m {F}{\Phi}^{-1}\phi(\overline{X})}.
  \end{equation}
  \item[(4)]  The angular metric tensors $\widehat{\hbar}$ and $\hbar$ are related by
 \begin{eqnarray}
    \widehat{\hbar}(\overline{X},\overline{Y}) &=& (m+1)  {F}^{2m}{\Phi}^{-2m} \big{\{} (\hbar(\overline{X},\overline{Y}) + m \ell(\overline{X})\,\ell(\overline{Y}) ) \\
   &&\nonumber
    +m F{\Phi}^{-1} ({F}{\Phi}^{-1}\,\phi(\overline{X})\,\phi(\overline{Y})  -\phi(\overline{X})\,\ell(\overline{Y})-\phi(\overline{Y})\,\ell(\overline{X}))\big{\}}.
      \end{eqnarray}
      \item[(5)]The relationship between the Finsler metric tensors $g$ and $\widehat{g}$ is given by
   \small{   \begin{eqnarray*}\label{gg3}
    \widehat{g}(\overline{X},\overline{Y})&=&  (m+1) {F}^{2m}{\Phi}^{-2m} \big{\{ } g(\overline{X},\overline{Y})+m  (2 m+1) {F}^{2} {\Phi}^{-2} 
      \,\phi(\overline{X})\,\phi(\overline{Y})
    \nonumber\\
    && +2 m \ell(\overline{X})\,\ell(\overline{Y}) -  2m{F}{\Phi}^{-1}( \phi(\overline{X})\,\ell(\overline{Y})+\phi(\overline{Y})\,\ell(\overline{X})) \big{\}}.
    \end{eqnarray*}}
\item[(6)] The Cartan torsions $\mathbf{T}$ and  $\widehat{\mathbf{T}}$ are related by 
\begin{eqnarray*}\label{tt3}
    2\widehat{\mathbf{T}}(\overline{X},\overline{Y},\overline{Z})&=&  2(m+1) \Phi^{-2 m} F^{2 m}\mathbf{T}(\overline{X},\overline{Y},\overline{Z})
    \\
    &&+2 m (m+1) \Phi^{-2 m} F^{2 m-1}\,\,\set{\hbar(\overline{X},\overline{Z})\,\ell(\overline{Y})
    +\hbar(\overline{Y},\overline{Z})\,\ell(\overline{X})} \nonumber\\
    &&-2 m (m+1) \Phi^{-2 m-1} F^{2 m}\,\set{\hbar(\overline{X},\overline{Z})\,\phi(\overline{Y})
    +\hbar(\overline{Y},\overline{Z})\,\phi(\overline{X})} \\
    && + (m+1)\left({D^\circ_{\gamma \overline{Z}}} \Phi^{-2 m} F^{2 m}\right)g(\overline{X},\overline{Y})\\
    &&+m (2 m+1) \left({D^\circ_{\gamma \overline{Z}}}\,\Phi^{-2 (m+1)} F^{2 m+2}\right) \,\phi(\overline{X})\,\phi(\overline{Y})
     \nonumber\\
    &&-2 m (m+1)\left({D^\circ_{\gamma \overline{Z}}}\, \Phi^{-2 m-1} F^{2 m+1}\right)\,\set{\phi(\overline{X})\,\ell(\overline{Y})+\phi(\overline{Y})\,\ell(\overline{X})}\\
    &&+2 m (m+1)\left({D^\circ_{\gamma \overline{Z}}}\, \Phi^{-2 m} F^{2 m}\right)\,\ell(\overline{X})\,\ell(\overline{Y}).
        \end{eqnarray*}
\end{description}
\end{prop}
\begin{proof} Let $(M,F)$ be a Finsler manifold equipped with a $\pi$-concurrent vector field $\overline{\varphi}$.  Under the $\phi$-concurrent generalized Kropina change {(\ref{change})},  we deduce the following, for all $Z \in  \mathfrak{X}(TM),\, \overline{U},\overline{X},\overline{Y} \in  \mathfrak{X}(\pi)$: \\
\noindent \textbf{(1)}  As the horizontal map ${\widehat{\beta}}$  of $\widehat{F}$  can be written in terms of the the horizontal map $\beta$ of $F$ (in the form  $\widehat{\beta}=\beta + \gamma \overline{U}$, for some  $\overline{U}$) and we have $\rho \circ \gamma =0$ in addition to 
 ${{{D}}}^{\circ}_{\gamma \overline{X}}\overline{Y}=\rho[\gamma
\overline{X}, {\beta} \overline{Y}]$, from \eqref{Berwald},  therefore, we obtain
   $$ {{\widehat{D}}}^{\circ}_{\gamma \overline{X}}\overline{Y}=\rho[\gamma
\overline{X}, \widehat{\beta}\, \overline{Y}]=\rho[\gamma \overline{X}, {\beta} \overline{Y}]
+\rho[\gamma \overline{X}, \gamma \overline{U}]=\rho[\gamma \overline{X}, {\beta} \overline{Y}]={{{D}}}^{\circ}_{\gamma \overline{X}}\overline{Y}.$$
\begin{eqnarray*}
\textbf{ (2)}\, Clearly, \,\,   d\widehat{E}(Z)&=&\frac 1 2 d\widehat{F}^2(Z)=\widehat{F}\,d\widehat{F}(Z)\\
  &=&\Phi^{-m} F^{m+1} \set{(m+1) \Phi^{-m} F^m\,dF(Z)-m \Phi^{-m-1} F^{m+1}\, d\Phi(Z)}\\
 &=& \Phi^{-2 m} F^{2 m+1}\,\{(m+1) dF(Z)-m \Phi^{-1} F\, d\Phi(Z)\}.\end{eqnarray*}
\noindent \textbf{(3)} From  Lemma \ref{B} \textbf{(a)}, \textbf{(b)}, it follows that
   \begin{eqnarray*}
                \widehat{\ell}(\overline{X})&=& d_J \widehat{F}(\widehat{\beta} \overline{X})=d_J \widehat{F}({\beta} \overline{X})=\frac{\partial \widehat{F}}{\partial F}\, d_J F(\beta \overline{X})
                 +\frac{\partial \widehat{F}}{\partial \Phi} \, d_J \Phi(\beta \overline{X})   \\
                  &=& (m+1) \Phi^{-m} F^m\,\ell(\overline{X})
     -m \Phi^{-m-1} F^{m+1}\, \phi(\overline{X}).
               \end{eqnarray*}
    \noindent \textbf{(4)} Using items \textbf{(1)}, \textbf{(3)}  above, Lemma \ref{B} \textbf{(a)}, \textbf{(b)}, \textbf{(e)} and Definition \ref{concurrent},  we get         
\begin{eqnarray*}
 \widehat{\hbar}(\overline{X},\overline{Y})&=&\widehat{F}({{\widehat{D}}}^{\circ}_{\gamma \overline{X}} \,\widehat{\ell})(\overline{Y})=
 \widehat{F}({{D}}^{\circ}_{\gamma \overline{X}} \,\widehat{\ell})(\overline{Y}) \\
   &=& \widehat{F}\, {{D}}^{\circ}_{\gamma \overline{X}} \,\set{(m+1) \Phi^{-m} F^m\,\ell(\overline{Y})
     -m \Phi^{-m-1} F^{m+1}\, \phi(\overline{Y})}\\
   &=& \widehat{F}\, \set{({{D}}^{\circ}_{\gamma \overline{X}} \,\left( (m+1) \Phi^{-m} F^m\right))\,\ell(\overline{Y})-({{D}}^{\circ}_{\gamma \overline{X}} \, m \Phi^{-m-1} F^{m+1})\, \phi(\overline{Y})}\\
   && +\widehat{F}\,\set{(m+1) \Phi^{-m} F^m\,({{D}}^{\circ}_{\gamma \overline{X}} \,\ell)(\overline{Y})-m \Phi^{-m-1}F^{m+1}
   ({{D}}^{\circ}_{\gamma \overline{X}} \, \phi)(\overline{Y})}\\
   &=&m (m+1) \, \Phi^{-2m} F^{2m+1}\, \{ (F^{-1}\,\ell(\overline{X})- \Phi^{-1} \, \mathbf{{B}}(\overline{X}))\,\ell(\overline{Y})\\
   &&  + (-\ell(\overline{X})+ \Phi^{-1} F\, \phi(\overline{X}) )\, \Phi^{-1} \,\phi(\overline{Y}) \}+(m+1)   \,\Phi^{-2m} F^{2m}\,\hbar(\overline{X},\overline{Y}).
\end{eqnarray*}
\noindent \textbf{(5)} Using items \textbf{(3), (4)} above and  the definition $\widehat{\hbar}:=\widehat{g}-\widehat{\ell} \otimes \widehat{\ell}$,  we obtain
    \begin{eqnarray*}
    \widehat{g}(\overline{X},\overline{Y}) &=&(m+1) \Phi^{-2 m} F^{2 m}\hbar(\overline{X},\overline{Y})
    +m (m+1) \Phi^{-2 (m+1)} F^{2 m+2} \,\phi(\overline{X})\,\phi(\overline{Y})  \nonumber\\
   &&-m (m+1) \Phi^{-2 m-1} F^{2 m+1}\,\set{\phi(\overline{X})\,\ell(\overline{Y})+\phi(\overline{Y})\,\ell(\overline{X})}\nonumber\\
   &&+m (m+1) \Phi^{-2 m} F^{2 m} \,\ell(\overline{X})\,\ell(\overline{Y}) \\
    &&+ \set{(m+1) \Phi^{-m} F^m\,\ell(\overline{X})
     -m \Phi^{-m-1} F^{m+1}\, \phi(\overline{X})}\, \times
     \\
     &&  \quad \set{(m+1) \Phi^{-m} F^m\,\ell(\overline{Y})
     -m \Phi^{-m-1} F^{m+1}\, \phi(\overline{Y})} \\
     &=&   (m+1) \Phi^{-2 m} F^{2 m}g(\overline{X},\overline{Y})+m (2 m+1) \Phi^{-2 (m+1)} F^{2 m+2} \,\phi(\overline{X})\,\phi(\overline{Y})
     \nonumber\\
    &&-2 m (m+1){F}^{2m}{\Phi}^{-2 m} \set{ \frac{F}{\Phi}\set{\phi(\overline{X})\,\ell(\overline{Y})+\phi(\overline{Y})\,\ell(\overline{X})}-\ell(\overline{X})\,\ell(\overline{Y})}.
    \end{eqnarray*}
\item[\textbf{(6)}] It follows from items \textbf{(1)}, \textbf{(5)} above and second part of \eqref{metricity}.
\end{proof}
\begin{thm}\label{F hat nondeg}
Let $(M, F)$ be a Finsler manifold admitting a concurrent $\pi$-vector field $\overline{\varphi}$ with associated  $\pi$-form $\phi$. The function $\widehat{F}$ defined by \eqref{change}  is  a Finsler metric if and only if 
\begin{equation}
\label{Eq:g_hat_non_degenerate}
mF^2  ||\overline{\varphi}||_{g}-(m-1)\Phi^2 \neq0.
\end{equation}
In other words,  the Finsler metric tensor $\widehat{g}$ of $\widehat{F}$ is non-degenerate if and only if the function $( mF^2  ||\overline{\varphi}||_{g}-(m-1)\Phi^2 )$  does not  vanish identically.
\end{thm}
 \begin{proof}
 The metric $\widehat{g}$ associated with $\widehat{F}$ is non-degenerate if and only if    $$\widehat{g}(\overline{U}, \overline{V})=0\,\, \, \,\forall\, \overline{U}\in\mathfrak{X}(\pi ) \implies \overline{V}=0.$$
 Assume that  $\widehat{g}(\overline{U}, \overline{V})=0,\,\, \, \forall\, \overline{U}\in\mathfrak{X}(\pi)$.  Then, relation \eqref{gg3} gives rise to
\begin{eqnarray}\label{sunsb}
   0 &=&  (m+1) \Phi^{-2 m} F^{2 m}g(\overline{U},\overline{V})+m (2 m+1) \Phi^{-2 (m+1)} F^{2 m+2} \,\phi(\overline{U})\,\phi(\overline{V})
     \nonumber\\
    &&-2 m (m+1) \Phi^{-2 m-1} F^{2 m+1}\,\set{\phi(\overline{U})\,\ell(\overline{V})+\phi(\overline{V})\,\ell(\overline{U})}\\  \nonumber
    &&+2 m (m+1) \Phi^{-2 m} F^{2 m}\,\ell(\overline{U})\,\ell(\overline{V}).
    \end{eqnarray}
Setting one time $\overline{W}=\overline{\varphi}$ then another time  $\overline{W}=\overline{\eta}$  in \eqref{sunsb},  we get
\begin{equation}\label{12q}
 \mathcal{ A}_1 \,\ell(\overline{Y}) + \mathcal{B}_1\, \phi({\overline{Y}})=0, \qquad \mathcal{A}_2\, \ell(\overline{Y}) +  \mathcal{ B}_2\, \phi({\overline{Y}})=0,
\end{equation}
where $\mathcal{A}_1 := 2 m (m+1) \Phi^{-2 m-1} F^{2 m-1} (\Phi^2-F^2 ||\overline{\varphi}||_{g})$,
\begin{eqnarray*}
   \mathcal{B}_1 &:=&\Phi^{-2 (m+1)} F^{2 m} \left(F^2 m (2 m+1) ||\overline{\varphi}||_{g}-\Phi^2 \left(2 m^2+m-1\right)\right),\\ \mathcal{A}_2 &:=& (m+1) \Phi^{-2 m} F^{2 m+1},\,  \qquad \mathcal{ B}_2 :=-m \,\Phi^{-2 m-1} F^{2 m+2}.
\end{eqnarray*}
This system of equations (\ref{12q}) has a non-trivial solution if and only if
$$  (m+1) \Phi^{-4 m-2} F^{4 m+1} \left({mF^2 ||\overline{\varphi}||_{g}-(m-1)\Phi^2}\right)=0. $$
That is, as $F\neq0$ over $\T M$,
$$mF^2 ||\overline{\varphi}||_{g}- (m-1)\Phi^2=0. $$
Consequently,
$$\overline{Y}\neq0 \quad \Longleftrightarrow \quad mF^2  ||\overline{\varphi}||_{g}-(m-1)\Phi^2 =0.$$
Therefore,  $\overline{Y}=0$ if and only if  the Finsler structure $F$ and the $\pi$-form $\Phi$ satisfy the condition
$$mF^2  ||\overline{\varphi}||_{g}-(m-1)\Phi^2 \neq0.$$
This means that the  $\phi$-generalized  Kropina metric tensor $\widehat{g}$ is non-degenerate if and only if the condition \eqref{Eq:g_hat_non_degenerate} is satisfied. Hence, the proof is completed.
\end{proof}
     Form now on, we consider that the   $\phi$-generalized  Kropina metric   $\widehat{F}$ satisfies the condition \eqref{Eq:g_hat_non_degenerate} and $\Phi^{m} >0$ in order  to make $\widehat{F}$ Finsler metric.  In fact,  $\widehat{F}$ is  a conic pseudo-Finsler metric. Let us give an example of the $\phi$-generalized  Kropina change of a Finsler metric which admits a $\pi$-concurrent vector field. Its calculations can be done by hand or using the Maple Finsler package \cite{NF_Package}. 
\begin{example}
\em{
Let $M=\{x:=(x^1,x^2,x^3)\in\Real^3:\, x^1>0\}$ and consider the conic Finsler metric defined on the domain 
$\mathcal{D}= \{(x,y:=y^1,y^2,y^3)) \in  TM |\, y^3  \neq 0   \}$ by
\begin{equation}\label{Fexample}
F(x,y)=x^1\sqrt{\frac{(y^1)^2y^3+(y^2)^3}{y^3} }.\end{equation}
Clearly, the non-vanishing components $g_{ij}$ of the Finsler metric tensor  are 
$$g_{11}=1, \quad g_{22}=\frac{3 (x^1)^2y^2}{y^3}, \quad g_{23}=-\frac{3}{2}\frac{ (x^1)^2(y^2)^2}{(y^3)^2},\qquad g_{33}=\frac{ (x^1)^2(y^2)^3}{(y^3)^3}.$$
Therefore, non-vanishing  components $C_{ijk}$ of the Cartan tensor are the following
$$C_{222} = \frac{3}{2}\frac{\, (x^1)^2 }{y^3}, \quad C_{223}=-\frac{3}{2}\frac{\, (x^1)^2 y^2}{(y^3)^2},\quad C_{233}= \frac{3}{2}\frac{\, (x^1)^2 (y^2)^2}{(y^3)^3}, \quad  C_{333}= -\frac{3}{2}\frac{\, (x^1)^2 (y^2)^3}{(y^3)^4}.$$
 In addition, the non-vanishing components $g^{ij}$ of the inverse metric tensor  are 
$$
  g^{11} =1, \quad g^{22} = \frac{4}{3}\frac{ y^3}{ (x^1)^2 y^2} , \quad  g^{23} =  \frac{ 2 (y^3)^2}{ (x^1)^2 (y^2)^2}  \quad  g^{33}=\frac{4}{3}\frac{ y^3}{ (x^1)^2 y^2} ,  \quad  g^{33}= \frac{ 4 (y^3)^3}{ (x^1)^2 (y^2)^3}   .
$$
Now, consider  the vector field whose  components given by 
$${\overline{\varphi}}^1(x)=x^1, \quad  {\overline{\varphi}}^2(x)={\overline{\varphi}}^3(x)=0.$$
Since  we have ${\overline{\varphi}}^i \,C_{ijk}=0$ and one can show that ${\overline{\varphi}}^i_{\,\,|j} =\delta _{i}^{j}$, we deduce that the Finsler metric $F$, defined by \eqref{Fexample}, admits a concurrent $\pi$-vector field given by $\overline{\varphi}={\overline{\varphi}}^i\,\overline{\partial_i}$, where $\overline{\partial_i}$ are the basis of fibres of $\pi^{-1}(TM)$. Hence the corresponding function $\Phi$ becomes $\Phi(x,y)= x^1 y^1$.

Therefore, we have
\begin{eqnarray*}
  \widehat{F}(x,y)&=&\frac{F^{m+1}(x,y)}{\Phi^m(x,y)}= \frac{x^1} {(y^1)^m} \left( \frac{(y^1)^2 y^3+(y^2)^3}{y^3} \right)^\frac{m+1}{2}
\end{eqnarray*}
 which defined on the domain: $$\mathbf{D}=\mathcal{D} \cap  \{(x,y) \in  TM |\, (y^1)^m  > 0   \}.$$ 
 Assume that $ \widehat{F}(x,y)$ satisfies the condition \eqref{Eq:g_hat_non_degenerate}, hence it is  a $\phi$-concurrent generalized Kropina  metric.   }
\end{example}
\begin{lem}\label{ddjEhat} Let $(M,F)$ be a Finsler manifold admitting concurrent $\pi$-vector field $\overline{\varphi}$.  Under the $\phi$-concurrent generalized Kropina  change {\em (\ref{change})}, we have: 
\begin{eqnarray}
   i_{\widehat{G}}\,dd_{J}\widehat{E}(X)	&=& \Phi^{-2 m-1} F^{2 m+1}\big{\{} (m+1) [   \Phi F^{-1} g(\overline{\mu},\rho X) - \Phi dF(X)\nonumber \\
    && +2 m \ell(\rho {X}) (F^{2} +\Phi^{-1}  F \ell(\overline{\mu}) -\phi(\overline{\mu})) -2 m \phi(\rho {X}) \ell(\overline{\mu}) ] \\ \nonumber 
    &&+m(2m+1)  \phi(\rho {X})[- \Phi^{-1}  F^3 + \Phi^{-1}  F \phi(\overline{\mu})]  +m F d\Phi (X)\big{\}}.
\end{eqnarray}
\end{lem}
\begin{proof}
 Because the difference between the two sprays constitutes a vertical vector field (i.e., $\widehat{G}=G+\gamma \overline{\mu}$, for some
$\pi$-vector field $\overline{\mu}$),  we find
\begin{equation}\label{ch5eq}
     \left.
    \begin{array}{rcl}
i_{\widehat{G}}\,(\frac{1}{2}\,dd_{J} \widehat{F}^2)(X)&=&i_{G+\gamma
\overline{\mu}}\,(\frac{1}{2}\,dd_{J} \widehat{F}^2)(X)=\frac{1}{2}i_{G}\,dd_{J}\widehat{F}^2(X)+\frac{1}{2}i_{\gamma \overline{\mu}}\,dd_{J}\widehat{F}^2(X).{\ \ \ \ }
 \end{array}
  \right.
 \end{equation}
Using Lemma \ref{B},  we derive
\small{\begin{eqnarray}\label{AB1}
              \frac{1}{2}\,i_{G}\,dd_{J}\widehat{F}^{2}(X) &=& \frac{1}{2} dd_{J}\widehat{F}^{2}(\beta \overline{\eta} ,X) \nonumber\\
 &=&\frac 1 2 \set{G \cdot d_J\widehat{F}^{2}(X) - X \cdot d_J\widehat{F}^{2}(G)-d_J\widehat{F}^{2}[G,X]} \nonumber\\
 &=&  G \cdot (\widehat{F} \widehat{\ell}(\rho X))-X \cdot( \widehat{F} \widehat{\ell}(\overline{\eta}))
 -\widehat{F}\widehat{\ell}(\rho[G,X]) \nonumber\\
 &=&{(G \cdot \widehat{F})\,\widehat{\ell}(\rho X)+\widehat{F}\, G \cdot \widehat{\ell}(\rho X)-(X \cdot \widehat{F}^2)
 -\widehat{F}\,\widehat{\ell}(\rho[G,X])}.
\end{eqnarray}}
Since, we have \begin{eqnarray*}
 G \cdot \widehat{F} &=& d\widehat{F}(G)= (m+1) \Phi^{-m} F^m\,dF(G)-m \Phi^{-m-1} F^{m+1}\, d\Phi(G) =m \Phi^{-m-1} F^{m+3}.\\
                                                     X \cdot \widehat{F}&=& d\widehat{F}(X)=\Phi^{-m} F^m \{ (m+1) \,dF(X)-m \Phi^{-1} F\, d\Phi(X)\}. \end{eqnarray*}  
Considering the above relations together with Lemma \ref{B} and Proposition \ref{hh1}, expression \eqref{AB1} simplifies to
\begin{eqnarray}\label{eq:ddjEhat} 
             \frac{1}{2}\,i_{G}\,dd_{J}\widehat{F}^{2}(X)
 &=&m \Phi^{-2m-1}\,F^{2m+3} \left((m+1) \,\ell(\rho {X})-m \Phi^{-1}F\, \phi(\rho {X})\right)\nonumber\\
 &&+\Phi^{-m}F^{m+1}\, G \cdot\left(\Phi^{-m}F^{m}\{(m+1)\,\ell(\rho {X})-m \Phi^{-1} F\, \phi(\rho {X})\}\right)\nonumber\\
 && -2\Phi^{-2m} F^{2m+1}\,\left((m+1) dF({X})-m \Phi^{-1} F\, d \Phi({X})\right)\nonumber\\
  &&-\Phi^{-2m} F^{2m+1}\,\, \left((m+1) \ell(\rho [G,X])-m \Phi^{-1} F  \phi(\rho [G,X])\right)\nonumber\\
 &=&\Phi^{-2 m-1} F^{2 m+1}\big{\{} (m+1)\left( 2 m F^{2} \ell(\rho {X})    - \Phi dF(X) \right)\nonumber\\&&-m (2 m+1) \Phi^{-1} F^{3} \phi(\rho {X})
 +m  F\, d\Phi(X) \big{\}}.
\end{eqnarray}
On the other hand, from \eqref{g}, we get
\begin{eqnarray}\label{AB3}
\frac{1}{2}i_{\gamma \overline{\mu}}\,dd_{J}\widehat{F}^2(X)&=&  \widehat{g}(\overline{\mu},\rho X) \nonumber\\
 &=& \Phi^{-2 m} F^{2 m}\big{\{} (m+1) g(\overline{\mu},\rho X)+m (2 m+1) \Phi^{-2} F^{2} \,\phi(\overline{\mu})\,\phi(\rho X)
     \nonumber\\
    &&-2 m (m+1) \Phi^{-1} F\,\set{\phi(\overline{\mu})\,\ell(\rho X)+\phi(\rho X)\,\ell(\overline{\mu})}\\ \nonumber 
    &&+2 m (m+1)\,\ell(\overline{\mu})\,\ell(\rho X)\big{\}} .
\end{eqnarray}
Substituting \eqref{AB3} and Formula \eqref{eq:ddjEhat}  into Equation (\ref{ch5eq}), it is evident that, following further calculations, \,
\begin{eqnarray*}\label{ABC}
     i_{\widehat{G}}\,dd_{J}\widehat{E}(X)&=&\Phi^{-2 m-1} F^{2 m+1}\big{\{} (m+1)\left( 2 m F^{2} \ell(\rho {X})    - \Phi dF(X) \right)\nonumber\\&&-m (2 m+1) \Phi^{-1} F^{3} \phi(\rho {X})
 +m  F\, d\Phi(X) \big{\}} \nonumber \\
&&+\Phi^{-2 m} F^{2 m}\big{\{} (m+1) g(\overline{\mu},\rho X)+m (2 m+1) \Phi^{-2} F^{2} \,\phi(\overline{\mu})\,\phi(\rho X)
     \nonumber\\ \nonumber
    &&-2 m (m+1) \Phi^{-1} F\,\set{\phi(\overline{\mu})\,\ell(\rho X)+\phi(\rho X)\,\ell(\overline{\mu})}+2 m (m+1)\,\ell(\overline{\mu})\,\ell(\rho X)\big{\}} 
    \\ \nonumber 
    &=&\Phi^{-2 m-1} F^{2 m+1}\big{\{} (m+1) [   \Phi F^{-1} g(\overline{\mu},\rho X) - \Phi dF(X)\nonumber \\
    && +2 m \ell(\rho {X}) (F^{2} +\Phi^{-1}  F \ell(\overline{\mu}) -\phi(\overline{\mu})) -2 m \phi(\rho {X}) \ell(\overline{\mu}) ] \\ \nonumber 
    &&+m(2m+1)  \phi(\rho {X})[- \Phi^{-1}  F^3 + \Phi^{-1}  F \phi(\overline{\mu})]  +m F d\Phi (X)\big{\}}.
\end{eqnarray*}
\vspace*{-1.2 cm}\[\qedhere\]
\end{proof}
\begin{thm}\label{th.22} 
Consider a Finsler manifold $(M,F)$ with a concurrent $\pi$-vector field $\overline{\varphi}$ and associated $\pi$-form $\phi$.   If $G$ is the geodesic spray of $F$,  then the geodesic spray $\widehat{G}$ of the $\phi$-concurrent generalized Kropina  metric $\widehat{F}$  is given by
\begin{equation}\label{geodesic sprays change}
\widehat{G}=G-\Psi_1 \,{\mathcal{C}}+\Psi_2 \,\gamma \overline{\varphi},
\end{equation}
where $\Psi_1:=\frac{2 m\Phi F^2 }{mF^2  ||\overline{\varphi}||_{g}-(m-1)\Phi^2 }$,
$\Psi_2:=\frac{mF^4 }{mF^2  ||\overline{\varphi}||_{g}-(m-1)\Phi^2}$.
\end{thm}
\begin{proof}
Given that the geodesic spray $\widehat{G}$ of the Finsler metric $\widehat{F}$ adheres to the equation \cite{r21} $$-d\widehat{E}= i_{\widehat{G}} dd_{J} \widehat{E}.$$ 
The expressions of $d\widehat{E}$ and $i_{\widehat{G}} dd_{J} \widehat{E}$ which are calculated in \eqref{AB} and Lemma \ref{ddjEhat}, respectively, leads to
 \small{
\begin{eqnarray*}- \set{(m+1)\,\Phi\,dF(X)-m \,F\, d\Phi (X)} &=&m F d\Phi (X)+(m+1) \big{\{ } \Phi F^{-1} g(\overline{\mu},\rho X) - \Phi dF(X)\nonumber \\
    && +2 m \ell(\rho {X}) (F^{2} +\Phi^{-1}  F \ell(\overline{\mu}) -\phi(\overline{\mu}))  \big{\}} \\ \nonumber 
    &&-2 m(m+1) \phi(\rho {X}) \ell(\overline{\mu})  \\ \nonumber 
    &&+m(2m+1)  \phi(\rho {X})[- \Phi^{-1}  F^3 + \Phi^{-1}  F \phi(\overline{\mu})] .
\end{eqnarray*}}
Which can be simplified to
\begin{eqnarray*}
 0&=& (m+1) [   \Phi F^{-1} g(\overline{\mu},\rho X)+2 m \ell(\rho {X}) (F^{2} +\Phi^{-1}  F \ell(\overline{\mu}) -\phi(\overline{\mu})) -2 m \phi(\rho {X}) \ell(\overline{\mu}) ] \\ \nonumber 
    &&+m(2m+1)  \phi(\rho {X})[- \Phi^{-1}  F^3 + \Phi^{-1}  F \phi(\overline{\mu})]  .
    \end{eqnarray*}
With the aid of the non-degenerate property of the Finsler metric $g$, the  above relationship gets shortened to
\begin{eqnarray}\label{ch52.eq.5}
   (m+1) \overline{\mu}&=&
   \{-2 m (m+1) \Phi^{-1} F^{2} +2 m (m+1) \Phi^{-1}  \phi({\overline{\mu}})-2 m (m+1) 
    F^{-1}\ell(\overline{\mu})\}\overline{\eta} \nonumber \\
   && +   \{m (2 m+1) \Phi^{-2 } F^{4}-m (2 m+1) \Phi^{-2} F^{2}\,\phi(\overline{\mu})
  \nonumber \\
   && +2 m (m+1) \Phi^{-1} F\,\ell(\overline{\mu})\}\overline{\varphi}.\,
\end{eqnarray}
Therefore,  the geometric objects $\ell(\overline{\mu})$ and $\phi(\overline{\mu})$ can be  determined by 
\begin{eqnarray}\label{mm}
    A_1 \,\ell(\overline{\mu})+B_1 \,\phi(\overline{\mu})&=&\,C_1 , \qquad
    A_2 \,\ell(\overline{\mu})+B_2 \:\phi(\overline{\mu})=\,C_2,
\end{eqnarray}
where
\begin{eqnarray*}
   A_1 &:=&  (m+1), \,\,\quad
  B_1 := - m \Phi^{-1} F,\,\,\quad C_1:= -m \Phi^{-1} F^{3},  \\
   A_2&:=& 2 m (m+1) \Phi^{-1} F^{-1} \left(\Phi^2-F^2 ||\overline{\varphi}||_{g}\right),\\
  B_2 &:=&\Phi^{-2}  \left(F^2 m (2 m+1) ||\overline{\varphi}||_{g}-\Phi^2 \left(2 m^2+m-1\right)\right),\\
  C_2&:=& m \Phi^{-2 } F^{2} \left(F^2 (2 m+1) ||\overline{\varphi}||_{g}-2 \Phi^2 (m+1)\right).  
\end{eqnarray*}
Applying the condition \eqref{Eq:g_hat_non_degenerate}, the  Algebraic system \eqref{mm} has  the solution 
\begin{eqnarray*}
  \ell(\overline{\mu}) &=& \frac{m\Phi F^3 }{(m-1)\Phi^2 -mF^2  ||\overline{\varphi}||_{g}}, \quad
  \phi(\overline{\mu}) =\frac{mF^2  \left(F^2 ||\overline{\varphi}||_{g}-2 \Phi^2\right)}{mF^2  ||\overline{\varphi}||_{g}-(m-1)\Phi^2 }.
\end{eqnarray*}
In light of Equation (\ref{ch52.eq.5}), which considers the assumption that $\widehat{G}=G+\gamma \overline{\mu}$, the canonical spray $\widehat{G}$ is provided by
\begin{eqnarray*}
\widehat{G}=G-\frac{2m \Phi F^2 }{mF^2  ||\overline{\varphi}||_{g}-(m-1)\Phi^2 } \,{\mathcal{C}}+\frac{mF^4 }{mF^2  ||\overline{\varphi}||_{g}-(m-1)\Phi^2 }\,\gamma \overline{\varphi}.
\end{eqnarray*}
Hence, the proof is completed.
\end{proof}
From \eqref{geodesic sprays change}, we note that under the $\phi$-concurrent generalized Kropina  change (\ref{change}) with $m\neq 0$,   the geodesic spray $G$ can not be invariant. In other words,  $\widehat{G} = G$ if and only if $\Psi _{1} =0=\Psi_{2}$ which means $ m\Phi F^2 = 0$ and $mF^4 =0$ which is impossible.
\begin{cor}\label{not projectively}
Let $(M,F)$ be a Finsler manifold admitting concurrent $\pi$-vector field $\overline{\varphi}$.   Under the $\phi$-concurrent generalized Kropina  change {\em (\ref{change})} with $m\neq 0$,   the geodesic sprays $G$ and  $\widehat{G}$ can never be projectively related.
\end{cor}
\begin{proof} Let $ \widehat{F}$ be the $\phi$-concurrent generalized Kropina  change of a Finsler metric $F$.  Assume that  $G$ and  $\widehat{G}$ are the geodesic sprays of $F$ and  $ \widehat{F}$, respectively.
Now, $G$ and $ \widehat{G}$ are projectively related if and only if $\widehat{G} = G - 2 \mathcal{P} \,\mathcal{C},$ where  $ \mathcal{P} $ is the projective factor which is positively homogeneous function of degree $1$ in  $y$.  Since $F $ is a non-zero function and $m\neq 0$, in view of the relation \eqref{geodesic sprays change}, we deduce that  $G$ and  $\widehat{G}$ are projectively related if and only if $\gamma \overline{\varphi}=0.$ Thus, 
$\overline{\varphi}=0$.  This contradicts our assumption that the $\pi$-vector field $\overline{\varphi}$ is everywhere nonzero.
\end{proof}
\begin{thm}\label{th.barthel} 
Let $(M,F)$ be a Finsler manifold admitting concurrent $\pi$-vector field $\overline{\varphi}$.  Under the $\phi$-concurrent generalized Kropina  change {\em (\ref{change})}, we have:  \,
\begin{description}
\item[(1)] 
The Barthel connections $\widehat{\Gamma}$ and $\Gamma$  are related by
\begin{equation}\label{ch52.eq.4}
         \widehat{\Gamma} =\Gamma +   \mathbb{F}, \qquad \mathbb{F}:=- \Psi_1\,J-d_J \Psi_1 \otimes \gamma \overline{\eta} - d_J \Psi_2 \otimes \gamma \overline{\varphi}.
\end{equation}
\item[(2)] The horizontal  projections $\widehat{h},\, h$ and vertical projections $\widehat{v},\,v$  are related, receptively, by
 $\widehat{h}=h+\frac 1 2 \mathbb{F}, \quad  \widehat{v}=v-\frac 1 2\mathbb{F}.$
\item[(3)] The Barthel curvature tensors   $\widehat{\Re}$ and $\Re$ are determined by
 $ \widehat{\Re} =\Re -\frac 1 2 [h,\mathbb{F}]-\frac 1 4 N_\mathbb{F}.$
\item[(4)] The horizontal counterpart of Berwald connection are related by 
  \begin{eqnarray*}
 \widehat{D^\circ}_{\widehat{\beta} \overline{X}} \, {\overline{Y}}&=&   {D^\circ}_{{\beta} \overline{X}} \overline{Y}-
 \frac 1 2\{\Psi_1\,D^\circ_{\gamma \overline{X}}\,\overline{Y}+d_J \Psi_1 ({\beta} \overline{X})\,
    D^\circ_{\gamma \overline{\eta}}\, \overline{Y}\\
    &&-d_J \Psi_1 ({\beta} \overline{X})\, \overline{Y}- d_J \Psi_1 (\beta \overline{Y})\,\overline{X}
    - d_J \Psi_2 ({\beta} \overline{X})\,D^\circ_{\gamma \overline{\varphi}}\, \overline{Y} \} \\
    && +\frac 1 2\set{dd_J \Psi_1 (\gamma \overline{Y},{\beta} \overline{X})\,  \overline{\eta} - dd_J \Psi_2 (\gamma \overline{Y}, {\beta} \overline{X})) \, \overline{\varphi}}.
\end{eqnarray*}
\end{description}
 \end{thm}
\begin{proof}
\begin{description}
\item[(1)]Based on Equation \eqref{fJformula}  and Formula \eqref{geodesic sprays change}, we can deduce 
\begin{eqnarray*}
  \widehat{\Gamma} &=& [J,\widehat{G}]
  = \left[J, G-\Psi_1\, \gamma \overline{\eta}+\Psi_2 \gamma \overline{\varphi} \right]
   =[J,G]+[\Psi_1\, \gamma \overline{\eta}-\Psi_2 \gamma \overline{\varphi},J]\\
  &=&  [J,G]+\Psi_1[\gamma \overline{\eta},J]+d\Psi_1\wedge i_{\gamma\overline{\eta}}\,J-d_{J}\Psi_1  \otimes \gamma \overline{\eta}\\
  && -\Psi_2[\gamma \overline{\varphi},J]-d\Psi_2\wedge i_{\gamma\overline{\varphi}}\,J+d_{J}\Psi_2 \otimes \gamma \overline{\varphi}.
\end{eqnarray*}
From \eqref{J properties}  and \eqref{P properties}, we deduce 
$$\widehat{\Gamma} =\Gamma -\Psi_1\,J-d_J \Psi_1 \otimes \gamma \overline{\eta} + d_J \Psi_2 \otimes \gamma \overline{\varphi}= \Gamma + \mathbb{F}.$$
\item[(2)] Applying  \eqref{hor. and ver. proj.} to the previous item, we get
 \begin{equation}\label{hor change}
\widehat{h} = \frac 1 2 (I +\widehat{\Gamma}) =  \frac 1 2 (I + \Gamma +\mathbb{F}) =h +\frac 1 2  \mathbb{F} \end{equation}
 and
$$\widehat{v}= \frac 1 2 (I -\widehat{\Gamma})= \frac 1 2 (I - \Gamma -\mathbb{F})=v- \frac 1 2 \mathbb{F}.$$
\item[(3)] The Barthel curvature tensor of $\widehat{F}$ is defined by $\widehat{\Re}=-\frac 1 2 [\widehat{h},\widehat{h}]$.  Now, item \textbf{(1)} above, \eqref{hor change} and formula
  \eqref{Nk} together with the aspects of the Fr\"{o}licher-Nijenhuis bracket lead to  \begin{eqnarray*}
 \widehat{\Re}&=&-\frac 1 2\, [h+\frac 1 2 \mathbb{F},h+\frac 1 2 \mathbb{F}] = -\frac{1}{2} \left([h,h] + \frac 1 2 [h,\mathbb{F}] +\frac 1 2 [\mathbb{F},h] +\frac 1 4 [\mathbb{F},\mathbb{F}]\right) \\
 &=& \Re -\frac 1 2 [h,\mathbb{F}]-\frac 1 4 N_\mathbb{F}.
  \end{eqnarray*}
\item[(4)]
Since  $v:=\gamma \circ K$, $h:=\beta \circ \rho$ and the Berwlad v-curvature $\widehat{S}^{\circ}=0$ together with Formulae \eqref{hor. and ver. proj.},  \eqref{brackect of two v.forms} and \eqref{JJ}, we obtain
 \begin{eqnarray*}
   \gamma \widehat{D^\circ}_{h {W}} \, \rho {Z}&=&\widehat{v}\,[ \widehat{h}W,JZ]\overset{\textbf{(2)}}{=}(v-\frac1 2 \mathbb{L})[(h+\frac 1 2\mathbb{L})W,JZ] \\
    &=& v[hW,JZ]+\frac 1 2 v[\mathbb{L}W,JZ]-\frac 1 2 \mathbb{L}[hW,JZ]-\frac 1 4 \mathbb{L}[\mathbb{L}W,JZ] \\
    &\overset{\eqref{Berwald}}{=}&  \gamma {D^\circ}_{h {W}} \overline{Z} +\frac \gamma 2 \big{\{}-\Psi_1\,K[\,JW,JZ]-d_J \Psi_1 (W)\, K[\, \gamma \overline{\eta},JZ] \\
    && + d_J \Psi_2 (W)\, K[\,\gamma \overline{\varphi},JZ] +(JZ \cdot \Psi_1)\,\rho W+(JZ \cdot d_J \Psi_1 (W))\,  \overline{\eta} 
    \\ &&- (JZ \cdot d_J \Psi_2 (W)) \, \overline{\varphi}
    +\Psi_1\,\rho([hW,JZ])+d_J \Psi_1 ([hW,JZ])  \, \overline{\eta} \\ &&
    - d_J \Psi_2 ([hW,JZ])  \, \overline{\varphi}\big{\}}\\
    &=&  \gamma {D^\circ}_{h {W}} \rho {Z}-\frac \gamma 2\{\Psi_1\,D^\circ_{J W}\,\rho Z+d_J \Psi_1 (W)\,
    D^\circ_{\gamma \overline{\eta}}\,\rho Z-d_J \Psi_1 (W)\,\rho Z\\
    &&
    - d_J \Psi_1 (Z)\,\rho W
    - d_J \Psi_2 (W)\,D^\circ_{\gamma \overline{\varphi}}\, \rho Z  +dd_J \Psi_1 (JZ,W)\,  \overline{\eta} \\
    &&- dd_J \Psi_2 (JZ, W) \, \overline{\varphi} \}.
    \end{eqnarray*}
Therefore,
\begin{eqnarray*}
 \widehat{D^\circ}_{\widehat{\beta} \,\overline{W}} \, {\overline{Z}}&=&   {D^\circ}_{{\beta} \overline{W}} \overline{Z}-
 \frac 1 2\{\Psi_1\,D^\circ_{\gamma \overline{W}}\,\overline{Z}+d_J \Psi_1 ({\beta} \overline{W})\,
    D^\circ_{\gamma \overline{\eta}}\, \overline{Z}\\
    &&-d_J \Psi_1 ({\beta} \overline{W})\, \overline{Z}- d_J \Psi_1 (\beta \overline{Z})\,\overline{W}
    - d_J \Psi_2 ({\beta} \overline{W})\,D^\circ_{\gamma \overline{\varphi}}\, \overline{Z}  \\
    && +dd_J \Psi_1 (\gamma \overline{Z},{\beta} \overline{W})\,  \overline{\eta} - dd_J \Psi_2 (\gamma \overline{Z}, {\beta} \overline{W}) \, \overline{\varphi}\}.
\end{eqnarray*}
\end{description}
\vspace*{-1.2 cm}\[\qedhere\]
\end{proof}
 \begin{cor}\label{cor:concurrent}
 Let $(M,F)$ be a Finsler manifold admitting a concurrent $\pi$-vector field $\overline{\varphi}$.  Under the  $\phi$-concurrent generalized Kropina  change \eqref{change}, we have:
 \begin{description}
 \item[(a)] The $\pi$-vector field $\overline{\varphi}$ is \textbf{not} concurrent with respect to the Finsler metric $\widehat{F}$.
 \item[(b)] If $\mathbb{F}=0$, then $\overline{\varphi}$ is concurrent with respect to the Finsler metric $\widehat{F}$.
\end{description}  
\end{cor}
\begin{proof}
  \begin{description}
 \item[(a)] Follows from Theorem \ref{th.barthel} \textbf{(4)} as $\widehat{D^\circ}_{\widehat{\beta} \,\overline{W}} \, {\overline{Z}} \neq {D^\circ}_{{\beta} \overline{W}} \overline{Z}$ together with Definition \ref{concurrent}.
 \item[(b)] Assume that $\mathbb{F}=0$,  then from Theorem \ref{th.barthel} \textbf{(1)}, we get $\widehat{\Gamma}= \Gamma$ which gives $\widehat{h}= h$ and  $\widehat{v}= v$. Thus, $ \widehat{D^\circ}_{\widehat{\beta} \,\overline{W}} \, {\overline{Z}} =  {D^\circ}_{\beta \,\overline{W}} \, {\overline{Z}}$. Moreover,  $ \widehat{D^\circ}_{\widehat{\gamma} \,\overline{W}} \, {\overline{Z}} =  {D^\circ}_{\gamma \,\overline{W}} \, {\overline{Z}}$ (by Proposition \ref{hh1} \textbf{(1)}). We deduce that if $\overline{\varphi}$ is concurrent with respect to $F$, then it is concurrent  with respect to  $\widehat{F}$.
 \end{description}  
 \vspace*{-1 cm}\[\qedhere\]
\end{proof}
\section{Almost rationally of  a $\phi$-Kropina change}
We end this paper by  a closer look of the effect of the $\phi$-concurrent generalized Kropina transformation on an almost rational Finsler metric. 
\begin{defn}\emph{ \cite{arFinslerTaha}}
A Finsler metric $F$ on  $M$ is said to be an almost rational Finsler metric if  its Finsler metric tensor  $g_{ij}(x,y)$ can be  written as a product of a positive smooth function $\theta$ on $TM$  and a symmetric non-degenerate  matrix $( a_{ij}(x,y))_{1 \leq i,j \leq n}$ with each of it entities $ a_{ij}(x,y)$ be a rational  in the directional argument $y$. That is,  $g_{ij}(x,y)= \theta (x,y)\, a_{ij}(x,y)$.
\par Additionally,   when $\theta $ is a rational function in $y$,   the Finsler metric $F$ is called a rational Finsler metric.
\end{defn}
\begin{rem}
The subsequent results warrant investigation based on the following facts:
\begin{description}
\item[(a)] A Riemannian metric is a quadratic Finsler metric; therefore, it is classified as a rational Finsler metric.  
\item[(b)]  The generalized Kropina metric, which can be viewed as a generalized Kropina change of a Riemannian metric, is listed as an almost rational Finsler metric, see \emph{\cite{arFinslerTaha}.}
\item[(c)] The $m$-th root metric falls into the category of  almost rational Finsler metrics \emph{\cite{arFinslerTaha}}. Moreover, the generalized Kropina change of an $m$-th root metric is an almost rational Finsler metric \emph{\cite{arFinslerTaha}.}
\end{description}
\end{rem}
\begin{thm}\label{L rational Finsler metric}
Given a rational Finsler metric  $F$ which admits a concurrent $\pi$-vector field $\overline{\varphi}$. Then the  $\phi$-concurrent generalized Kropina  metric $\widehat{F}$ is \\
-- a rational Finsler metric provided that $m \in \mathbb{Z}$.\\
-- an almost rational Finsler metric provided that $m$ not an integer.
\end{thm}
\begin{proof}
Since $F$ is a rational Finsler metric, all its Finsler metric tensor components $g_{ij}(x,y)$ are $y$-rational functions and can be expressed in the form $g_{ij}(x,y)=\zeta(x,y)\, a_{ij}(x,y)$ such that both $\zeta$ and $a_{ij}$ are $y$-rational functions. Also, $F^2 = g_{ij}y^i y^j$ is $y$-rational function. Furthermore, $F \ell_{i}= g_{ri} y^r $ are $y$-rational functions $\forall i=1,...,n$. 
Similarly,  the functions $$\ell_{i}\,\ell_{j} = \frac{g_{ri} y^r}{F} \frac{g_{kj} y^k}{F} =  \frac{g_{ri} y^r g_{kj} y^k}{F^2} = \frac{\zeta a_{ri} y^r  \zeta a_{kj} y^k}{\zeta a_{ms}y^m y^s} =\zeta \frac{ a_{ri} y^r  a_{kj} y^k}{a_{ms}y^m y^s} $$  are always $y$-rational. \\Now the local expression of the $\phi$-concurrent generalized Kropina  metric tensor, which can be deduced from the formula \eqref{gg3}, is given by
\begin{eqnarray*}\label{gij local}
    \widehat{g}_{ij}&=&  (m+1) \Phi^{-2 m} F^{2 m}g_{ij}+m (2 m+1) \Phi^{-2 (m+1)} F^{2 m+2} \,\phi_{i}\,\phi_{j}
     \nonumber\\
    &&-2 m (m+1) \Phi^{-2 m-1} F^{2 m+1}\,\set{\phi_{i}\,\ell_{j}+\phi_{j}\,\ell_{i}}+2 m (m+1) \Phi^{-2 m} F^{2 m}\,\ell_{i}\,\ell_{j}
    \nonumber\\
    &=&m(m+1)\left(\frac{F}{\Phi}\right)^{2m} \set{ \frac{1}{m}g_{ij}+ \frac{(2 m+1)}{m+1} \frac{F^{2 }}{\Phi^{2}} \,\phi_{i}\,\phi_{j} +2\ell_{i}\,\ell_{j} -\frac{2F}{\Phi}\set{\phi_{i}\,\ell_{j}+\phi_{j}\,\ell_{i}} }\\
    &=& \widehat{\zeta} \,\widehat{a}_{ij}
    \end{eqnarray*}
 provided  that  $ \widehat{\zeta} = m(m+1)\left(\frac{F^{2 }}{\Phi^{2}}\right)^{m}$ and $$\widehat{a}_{ij}=\frac{1}{m}g_{ij}+ \frac{(2 m+1)}{m+1} \frac{F^{2 }}{\Phi^{2}} \,\phi_{i}\,\phi_{j} +2\ell_{i}\,\ell_{j} -\frac{2F}{\Phi}\set{\phi_{i}\,\ell_{j}+\phi_{j}\,\ell_{i}} .$$ Thereby, when $m$ is an integer, the function $ \widehat{\zeta}$ is $y$-rational and for other values of $m$, the function $ \widehat{\zeta}$ is not $y$-rational. Moreover, based on the above discussion the functions $\widehat{a}_{ij}$ are always $y$-rational.
 Hence, the  $\phi$-concurrent generalized Kropina  metric $\widehat{F}$ is a rational Finsler metric provided that $m \in \mathbb{Z}$ and for other values of $m$, it is an almost rational Finsler metric.
\end{proof} 
\begin{thm}\label{L AR Finsler metric}
Given an almost rational Finsler metric  $F$ which admits a concurrent $\pi$-vector field $\overline{\varphi}$. Then the  $\phi$-concurrent generalized Kropina metric $\widehat{F}$ is  an almost rational Finsler metric.
\end{thm}
\begin{proof}
As $F$ is an almost rational Finsler metric implies,  by defintion,  $g_{ij} = \zeta \, a_{ij}$ with $\zeta$ be not $y$-rational function. The  $\phi$-concurrent generalized Kropina  metric tensor can be written in the form
\begin{eqnarray*}\label{ar gij local}
    \widehat{g}_{ij}
    &=&m(m+1)\,\left(\frac{F^{2 }}{\Phi^{2}}\right)^{m} \\&& \times
    \set{ \frac{1}{m}g_{ij}+ \frac{(2 m+1)}{m+1} \frac{F^{2 }}{\Phi^{2}} \,\phi_{i}\,\phi_{j} +2\ell_{i}\,\ell_{j} -\frac{2F}{\Phi}\set{\phi_{i}\,\ell_{j}+\phi_{j}\,\ell_{i}} }\\
    &=&
    m(m+1) \left(\frac{F^{2 }}{\Phi^{2}}\right)^{m}
    \\&& \times
   \set{ \frac{1}{m}\zeta \,a_{ij}+ \frac{(2 m+1)}{m+1} \frac{\zeta \, a_{rs}y^r y^s}{\Phi^{2}} \,\phi_{i}\,\phi_{j} +2\zeta \frac{ a_{ri} y^r  a_{kj} y^k}{a_{ms}y^m y^s}-\frac{2\zeta \, y^r}{\Phi}\set{\phi_{i}\,a_{rj}+\phi_{j}\,a_{ir}} }\\
    &=&
    m(m+1)\,\zeta \left(\frac{F^{2 }}{\Phi^{2}}\right)^{m}
    \\&& \times \set{\frac{1}{m}a_{ij}+ \frac{(2 m+1)}{m+1} \frac{a_{rs}y^r y^s}{\Phi^{2}} \,\phi_{i}\,\phi_{j} +2 \frac{ a_{ri} y^r  a_{kj} y^k}{a_{ms}y^m y^s}-\frac{2  y^r}{\Phi}\set{\phi_{i}\,a_{rj}+\phi_{j}\,a_{ir}} }.
    \end{eqnarray*}
    That is, $$ \widehat{g}_{ij} (x,y)=\widehat{\zeta}(x,y)\, \widehat{\textbf{a}}_{ij}(x,y)$$
where \begin{equation}\label{etaAR}
\widehat{\zeta}=m(m+1)\,\zeta \left(\frac{F^{2 }}{\Phi^{2}}\right)^{m}
\end{equation}
and
 \begin{equation}
\widehat{\textbf{a}}_{ij}=\frac{1}{m}a_{ij}+ \frac{(2 m+1)}{m+1} \frac{a_{rs}y^r y^s}{\Phi^{2}} \,\phi_{i}\,\phi_{j} +2 \frac{ a_{ri} y^r  a_{kj} y^k}{a_{ms}y^m y^s}-\frac{2  y^r}{\Phi}\set{\phi_{i}\,a_{rj}+\phi_{j}\,a_{ir}}.
\end{equation}
Since $\zeta$ is not $y$-rational function,  the function $ \widehat{\zeta}$, defined by \eqref{etaAR}, maybe or maybe not $y$-rational function.  However,  the functions $\widehat{\textbf{a}}_{ij}$ are always $y$-rational.
 Hence, the  $\phi$-concurrent generalized Kropina  metric $\widehat{F}$ is an almost rational Finsler metric.
\end{proof}
\section*{Conclusion}
We have investigate what we call the  $\phi$-concurrent generalized Kropina change $\widehat{F}$ of an arbitrary Finsler metric $F$ and find the associated Finslerian geometric objects of $\widehat{F}$ in terms of those of $F$ in Proposition \ref{hh1} and Theorems \ref{th.22},  \ref{th.barthel}.  Also, we give a necessary and sufficient condition to make $\widehat{F}$  Finsler metric (in Theorem \ref{F hat nondeg}).  This enable us to prove that:
\begin{description}
\item[(a)]The geodesic sprays $G$ and  $\widehat{G}$ can never be projectively related (see, Corollary \ref{not projectively}).
\item[(b)] If the $(1,1)$-tensor $\mathbb{F}$ defined in \eqref{ch52.eq.4} vanishes identically, the the Barthel connection, horizontal and vertical projectors, the curvature of  Barthel connection and Berwald connections are invariant. Therefore,  the vector filed $\overline{\phi}$ becomes   concurrent with respect to $\widehat{F}$ (see, Corollary \ref{cor:concurrent}).
\item[(c)] As a coordinate-study application (\S 3) of our results,  we deduce that the  $\phi$-concurrent generalized Kropina change preserves the almost rational property of the initial Finsler metric ${F}$ (in Theorem \ref{L AR Finsler metric}).
\item[(d)] Also, we have noted that from this paper and \cite{square metric, Soleiman-Taha_Mat} the following:\\
-- The vertical counterpart for Berwald connection is invariant \eqref{D}.\\
-- Given a Finsler manifold  $(M,F)$ equipped with a $\pi$-concurrent vector field $\overline{\varphi}$.  Under any change  $\widehat{F}  \longrightarrow F$ gives the geodesic spray $G$ change $\widehat{G}=G-\Psi_1 \,{\mathcal{C}}+\Psi_2 \,\gamma \overline{\varphi},$ where $\Psi_1 ,\Psi_1$ arbitrary smooth function on $TM$,  Theorem  \ref{th.barthel} and Corollary  \ref{cor:concurrent} still hold. 
\item[(e)] Finally,  we prove that the  $\phi$-concurrent generalized Kropina change ($F  \longrightarrow \widehat{F}$) preserves the almost rational property of the initial Finsler metric ${F}$ in \S 3.
\end{description}

\end{document}